\documentclass[12pt]{amsart}
\usepackage{etex}
\usepackage{geometry}
\usepackage{color}
\usepackage{todonotes}
\usepackage{amsfonts, amssymb, amscd}
\numberwithin{equation}{section}
\usepackage[symbol]{footmisc}

\usepackage{tikz}
\usetikzlibrary{matrix}
\usetikzlibrary{decorations.pathreplacing}
\usepackage{bm}
\usepackage{verbatim}
\usepackage{mathrsfs}
\usepackage{graphicx}
\usepackage{tikz-cd}
\usepackage{subcaption}
\usepackage{listings}
\usepackage{subfiles}
\usepackage[toc,page]{appendix}
\usepackage{mathtools}
\usepackage{comment}
\usepackage{enumerate}
\usepackage{enumitem}
\usepackage[all]{xy}
\usepackage{arydshln}
\graphicspath{{images/}}
\usepackage{appendix}
\usepackage{hyperref}
\hypersetup{
	colorlinks=true,
	citecolor=red,
	linkcolor=blue,
	filecolor=magenta,
	urlcolor=red,
}
\newtheorem{thm}{Theorem}[section]

\newtheorem{cor}[thm]{Corollary}
\newtheorem{lem}[thm]{Lemma}
\newtheorem{prop}[thm]{Proposition}

\theoremstyle{definition}

\theoremstyle{definition}

\begin{document}
	\title[Coincidence of the Bruhat order and the secondary Bruhat order]{The coincidence of the Bruhat order and the secondary Bruhat order on $\mathcal{A}(n,k)$}
	\author{Tao Zhang and Houyi Yu$^{\ast}$}  \thanks{*Corresponding author}
	\keywords{Bruhat order, secondary Bruhat order, $(0,1)$-matrix.}
	\begin{abstract}
		Given a positive integer $n$ and a nonnegative integer $k$ with $k\leq n$,
we denote by $\mathcal{A}(n,k)$ the class of all $n$-by-$n$ $(0,1)$-matrices with constant row and column sums $k$.
In this paper, we show that the Bruhat order and the secondary Bruhat order coincide on $\mathcal{A}(n,k)$ if and only if either $0\leq n\leq 5$ or $k\in\{0,1,2,n-2,n-1,n\}$ with $n\geq 6$.
	\end{abstract}
	
	\address{School of Mathematics and Statistics, Southwest University, Chongqing 400715, China}
	\email{2816187279@qq.com}
	
	\address{School of Mathematics and Statistics, Southwest University, Chongqing 400715, China}
	\email{yuhouyi@swu.edu.cn}
	\maketitle
	
\section{Introduction}
Let $m$ and $n$ be two positive integers and let $R=(r_1,r_2,\ldots,r_m)$ and $S=(s_1,s_2,\ldots,s_n)$ be two sequences of positive integers with the same sum.
Write $\mathcal{A}(R,S)$ for the class of all $m$-by-$n$ $(0,1)$-matrices with row sum vector $R$ and column sum vector $S$. For notational simplicity, we write simply $\mathcal{A}(n,k)$ for the class of all $n$-by-$n$ $(0,1)$-matrices with constant row and column sums $k$.

In recent years, there has been a considerable amount of interest in the combinatorial properties of the class $\mathcal{A}(R,S)$. For example, Gale \cite{Gal57} and Ryser \cite{Ry57} independently discovered a necessary and sufficient condition for $\mathcal{A}(R,S)$ to be nonempty.
Brualdi and Hwang \cite{BH04} extended the Bruhat order on the symmetric group $\mathfrak{S}_n$ to any nonempty classes $\mathcal{A}(R,S)$, still called the Bruhat order, which is a refinement of the secondary Bruhat order formally introduced a few years later by Brualdi and Deaett \cite{BD07}, building on the work of \cite{BH04}.
These two partial orders have been studied extensively in the literature  (see, e.g., \cite{Bru06,CFM13,CFM12,FCS20,Ghe14} and references therein).

In this work we classify all classes $\mathcal{A}(n,k)$ where the Bruhat order and the secondary Bruhat order coincide.
Before stating our main result, we recall a few definitions.

Let $\mathcal{A}(R,S)$ be a nonempty class of $(0,1)$-matrices.
Without loss of generality, we assume that $R$ and $S$ are nonincreasing.
If $t_1>t_2>\cdots>t_l$ are positive integers such that $\{t_1,t_2,\cdots,t_l\}=\{r_1,r_2,\cdots,r_m\}$, then we write $R=(t_1^{i_1},t_2^{i_2},\cdots,t_l^{i_l})$ where $i_k$ are the multiplicities of $t_k$ for $1\leq k\leq l$.
Given an $m$-by-$n$ matrix $A=[a_{ij}]$, let $\Sigma_A=[\sigma_{i j}(A)]$ denote the $m$-by-$n$ matrix where
  \begin{align*}
  	\sigma_{i j}(A)=\sum_{k=1}^{i} \sum_{l=1}^{j} a_{kl} \quad(1 \leqslant i \leqslant m, 1 \leqslant j \leqslant n).
  \end{align*}
The Bruhat order on $\mathcal{A}(R,S)$ is defined as follows.
For two $(0,1)$-matrices $A$ and $C$ in $\mathcal{A}(R,S)$, we say that $A$ precedes $C$ in the \emph{Bruhat order}, written $A\preceq_B C$, if $\Sigma_A\geq \Sigma_C$ by the entrywise order, that is, $\sigma_{ij}(A)\geq \sigma_{ij}(C)$ for all $i$ and $j$ with $1 \leqslant i \leqslant m$ and $1 \leqslant j \leqslant n$.
This is a generalization of the well-known Bruhat order for permutation matrices (equivalently, for the symmetric groups), see, e.g., \cite{BB05}.

The secondary Bruhat order on $\mathcal{A}(R,S)$ is based on the notion of interchanges.
Given an element $A$ in $\mathcal{A}(R,S)$,
an \emph{interchange} of $A$ is a sequence of replacements of a submatrix of type $L_2$ by type $I_2$ or vice versa, where
$$
  I_{2}=\begin{bmatrix}
  	1 & 0 \\
  	0 & 1
 \end{bmatrix}\quad \text{and}\quad L_{2}=\begin{bmatrix}
  	0 & 1 \\
  	1 & 0
 \end{bmatrix}.
 $$
Observe that an interchange always takes a matrix in a class $\mathcal{A}(R,S)$ to another matrix in the same class,
and that any two matrices in $\mathcal{A}(R,S)$ can be transformed into each other by a finite sequence of interchanges \cite{Bru06,Ry64}.
For two matrices $A$ and $C$ in $\mathcal{A}(R,S)$, we say $A$ precedes $C$ in the \emph{secondary Bruhat order}, written $A \preceq_{\widehat{B}} C$, provided that $C$ can be converted into $A$ by a sequence of one-sided interchanges $L_2 \rightarrow I_2 $ that replace a submatrix of order 2 equal to $L_2$ by $I_2$.

It is not difficult to see that the Bruhat order is a refinement of the secondary Bruhat order.
Brualdi and Deaett \cite{BD07} proved that the two Bruhat orders coincide on $\mathcal{A}(n,1)$ and $\mathcal{A}(n,2)$, but they are different on $\mathcal{A}(6,3)$, from which we see that these two orders are different in general.
Thus, it would be interesting to characterize all classes $\mathcal{A}(R,S)$ where the two Bruhat orders coincide.
Most recently, Fernandes, Cruz and Salom\~ao \cite{FCS20} showed that the two orders coincide on $\mathcal{A}(R,S)$ with $R=(1^m)$ or $R=(2^m)$, but they are different on $\mathcal{A}(n,3)$ for $n\geq6$.
According to \cite[Proposition 6]{FCS20}, we see that the two orders are the same on $\mathcal{A}(n,n-2)$ and $\mathcal{A}(n,n-1)$ for $n\geq6$.
Because both $\mathcal{A}(n,0)$ and $\mathcal{A}(n,n)$ consist of only one element, the two orders also coincide on these two classes of $(0,1)$-matrices.
Summing up these
conclusions yields that the Bruhat order and the secondary Bruhat order coincide on $\mathcal{A}(n,k)$ if either $0\leq n\leq 5$ or $k\in\{0,1,2,n-2,n-1,n\}$ with $n\geq 6$.
In this note we show that the converse is also true.

The paper is organized as follows. After summarizing some necessary background and basic facts on posets and $(0,1)$-matrices needed in the sequel, we present in Section \ref{sec:prelim} our main result, that is, Theorem \ref{thm:mainequal}, which provides a complete classification of the classes $\mathcal{A}(n,k)$ where  the Bruhat order agrees with the secondary Bruhat order. Section \ref{Sec:proofmainresult} is then devoted to the proof of this theorem.

	\section{Preliminaries}\label{sec:prelim}
We begin with reviewing the background information on posets and $(0,1)$-matrices. Details can be found for example in \cite{Sta12,Bru06}.
Unless otherwise specified, $n$ and $k$ are assumed to be nonnegative integers such that $k\leq n$.
Let $(P,\leq)$ be a finite partially ordered set (or poset, for short). For any $a, b \in P$, we use the notation $a < b$ to mean $a \leq b$, but $a \neq b$. Moreover, we say $b$ \emph{covers} $a$ or $a$ is \emph{covered} by $b$, denoted $a\lessdot b$ or $b\gtrdot a$, if $a<b$ and there exists no element $c \in P$ such that $a < c< b$.
A finite poset is determined uniquely by its cover relations.
Given an $n$-by-$n$ matrix $A$, we write $A[\{i_1, i_2,\cdots,i_r\}, \{j_1,j_2,\cdots,j_s\}]$ for the submatrix of $A$ that lies in rows $i_1, i_2,\cdots,i_r$ and columns $j_1,j_2,\cdots,j_s$.

In \cite{BD07}, Brualdi and  Deaett characterized the cover relations of the secondary Bruhat order on classes $\mathcal{A}(R,S)$.
	
	\begin{lem}\cite[Theorem 3.1]{BD07}\label{lem:coversecond} Let $A=\left[a_{i j}\right]$ be a matrix in $\mathcal{A}(R,S)$ where $A[\{i, j\}, \{k, l\}]=L_{2}$. Let $A^{\prime}=\left[a_{i j}^{\prime}\right]$ be the matrix obtained from $A$ by the $L_{2}$$ \rightarrow $$I_{2}$  interchange that replaces $A[\{i, j\}, \{k, l\}]=L_{2}$ with $I_{2}$. Then $A$ covers $A^{\prime}$ in the secondary Bruhat order on $\mathcal{A}(R, S)$ if and only if
\begin{enumerate}
\item[\rm (1)] $a_{p k}=a_{p l}$  $(i<p<j)$,
\item[{\rm (2)}] $a_{i q}=a_{j q}$  $(k<q<l)$,
\item[{\rm (3)}] $a_{p k}=0$ and $a_{i q}=0 $ imply $a_{p q}=0$ $(i<p<j, k<q<l)$, and
\item[{\rm (4)}]  $a_{p k}=1$ and $a_{i q}=1 $ imply $a_{p q}=1$ $(i<p<j, k<q<l)$.
\end{enumerate}
	\end{lem}

A \emph{partition} of a nonnegative integer $n$ is a sequence of nonnegative integers $\lambda=(\lambda_1,\lambda_2,\ldots,\lambda_k)$ satisfying
$\lambda_1\geq \lambda_2\geq\cdots\geq \lambda_k$ and $\sum_{i=1}^k\lambda_i=n$. Any $\lambda_i=0$ is considered irrelevant, and we identify $\lambda$ with the infinite sequence
$\lambda=(\lambda_1,\lambda_2,\ldots,\lambda_k,0,0,\ldots)$.
Let $\lambda$ be a partition of $n$. The \emph{conjugate partition} of $\lambda$ is the partition $\lambda^*=(\lambda_1',\lambda_2',\ldots)$, where
$\lambda_i'=|\{j:\lambda_j\geq i\}|$.
Let $\lambda$ and $\mu$ be two partitions of the same nonnegative integer. We say that $\lambda$ precedes $\mu$ in the \emph{dominance order}, denoted $\lambda \preceq \mu$,
if
$$
\lambda_1+\lambda_2+\cdots+\lambda_i\leq \mu_1+\mu_2+\cdots+\mu_i
$$
for all $i\geq1$.
Given a positive integral vector $R=(r_1,r_2,\ldots,r_m)$, we denote by $\widehat{R}$ the partition obtained from $R$ by reordering its coordinates decreasingly.

The well-known Gale-Ryser theorem, due independently to Gale \cite{Gal57} and Ryser \cite{Ry57}, provides a description for the existence
of $(0,1)$-matrices with prescribed row and column sums. See \cite{Bru06} and \cite{Ry64} for more details.

	\begin{lem}\label{lem:existm} (Gale-Ryser) Let $R=\left(r_{1}, r_{2},\ldots,r_{m}\right) $ and  $S=\left(s_{1}, s_{2},\ldots, s_{n}\right) $ be two nonnegative integral vectors
with the same sum.
Then $\mathcal{A} (R, S)$ is nonempty if and only if $\widehat{S} \prec (\widehat{R})^{*}$.
	\end{lem}

         In order to show the coincidence of the Bruhat order and the secondary Bruhat order on $\mathcal{A}(n,3)$ for $n\in\{3,4,5\}$, Fernandes, Cruz and Salom\~ao \cite{FCS20} presented an elegant duality property of $\mathcal{A}(R,S)$.
	\begin{lem}\cite[Proposition 6]{FCS20}\label{lem:dualthm}
Let $R=\left(r_{1}, \ldots, r_{m}\right) $, and  $S=\left(s_{1}, \ldots, s_{n}\right)$  be two nonincreasing positive integral vectors such that  $\mathcal{A}(R, S)\neq \emptyset$.  Let  $U=\left(n-r_{m}, \ldots, n-r_{1}\right)$, and  $Q=\left(m-s_{n}, \ldots, m-s_{1}\right)$. If the Bruhat order and the secondary Bruhat order coincide on $\mathcal{A}(R, S)$, then these two orders coincide on  $\mathcal{A}(U, Q)$.
	\end{lem}
       The following result is an immediate corollary of Lemma \ref{lem:dualthm}.
	\begin{cor}\label{lem:dualnkn-k}
		The Bruhat order and the secondary Bruhat order coincide on $\mathcal{A}(n,k)$ if and only if these two orders coincide on $\mathcal{A}(n,n-k)$.
	\end{cor}

The following result follows from \cite[Theorem 4.4]{BD07} and \cite[Section 3]{FCS20}, which together with Corollary \ref{lem:dualnkn-k}, indeed classifies all classes $\mathcal{A}(n,k)$ where the Bruhat order coincides with
the secondary Bruhat order, as stated in Theorem \ref{thm:mainequal}.

	\begin{lem}\label{lem:knocoincidenk012} The Bruhat order agrees with the secondary Bruhat order on $\mathcal{A}(n,0)$, $\mathcal{A}(n,1)$, $\mathcal{A}(n,2)$,
		$\mathcal{A}(3,3)$, $\mathcal{A}(4,3)$ and $\mathcal{A}(5,3)$. These two orders do not coincide on $\mathcal{A}(n,3)$, for $n \geq 6$.
	\end{lem}

\begin{cor}\label{cor:equal}
If $0\leq n\leq 5$, or $k\in\{0,1,2,n-2,n-1,n\}$ with $n\geq6$, then the Bruhat order coincides with the secondary Bruhat order on $\mathcal{A}(n,k)$.
\end{cor}
\begin{proof}
This result is a straightforward consequence of Corollary \ref{lem:dualnkn-k} and Lemma \ref{lem:knocoincidenk012}.
\end{proof}

We are now in a position to give our main result, whose proof needs more preparation and will be provided in Section \ref{Sec:proofmainresult}.

\begin{thm}\label{thm:mainequal}	
The Bruhat order and the secondary Bruhat order coincide on $\mathcal{A}(n,k)$ if and only if either $0\leq n\leq 5$ or $k\in\{0,1,2,n-2,n-1,n\}$ with $n\geq6$.
\end{thm}

	\section{Proof of Theorem \ref{thm:mainequal}}\label{Sec:proofmainresult}

Our arguments are inspired by an example, due to Brualdi and Deaett \cite{BD07}, showing that there exist in $\mathcal{A}(n,3)$ two different matrices that are incomparable in the secondary Bruhat order but one of which is less than the other in the Bruhat order, and thereby proving these two orders are not the same on $\mathcal{A}(n,3)$.
The key matrices used in \cite{BD07} are
\begin{align*}
				A=\begin{bmatrix}
					1 & 0 & 0 & 0\\
					1 & 0 & 1 & 1\\
					1 & 1 & 0 & 1\\
					0 & 0 & 0 & 1
				\end{bmatrix}, \quad C=\begin{bmatrix}
					0 & 0 & 0 & 1\\
					1 & 0 & 1 & 1\\
					1 & 1 & 0 & 1\\
					1 & 0 & 0 & 0
				\end{bmatrix}\quad\text{and}\quad
				\quad D=\begin{bmatrix}
					0 & 0 & 0 & 1\\
					1 & 1 & 0 & 1\\
					1 & 0 & 1 & 1\\
					1 & 0 & 0 & 0
				\end{bmatrix}.
			\end{align*}

\begin{lem}\label{lem:xvg123}
If there exist $(0,1)$-matrices $G_{1}$, $G_{2}$ and $G_{3}$  such that
				\begin{align*}
					V_n=\left[\begin{array}{c|c}
						V & G_{1}\\
						\hline G_{2}& G_{3}
					\end{array}\right]
				\end{align*}
are elements of $\mathcal{A}(n,k)$ for all $V \in\{A, C, D\}$,
then the Bruhat order and secondary Bruhat order do not coincide on $\mathcal{A}(n,k)$.
\end{lem}

\begin{proof} By calculation we get
\begin{align*}
				\Sigma_{A}=\begin{bmatrix}
					1 & 1 & 1 & 1\\
					2 & 2 & 3 & 4\\
					3 & 4 & 5 & 7\\
                    3 & 4 & 5 & 8
				\end{bmatrix}, \quad
				\Sigma_{C}=\begin{bmatrix}
					0 & 0 & 0 & 1\\
					1 & 1 & 2 & 4\\
					2 & 3 & 4 & 7\\
					3 & 4 & 5 & 8
					\end{bmatrix}
			\quad \text{and}\quad
				\Sigma_{D}=\begin{bmatrix}
					0 & 0 & 0 & 1\\
					1 & 2 & 2 & 4\\
					2 & 3 & 4 & 7\\
					3 & 4 & 5 & 8
					\end{bmatrix},
			\end{align*}
showing that $\Sigma_{A}> \Sigma_{D}> \Sigma_{C}$.
Since $G_1$, $G_2$ and $G_3$ are common parts of  $A_n$, $C_n$ and $D_n$, we have $\Sigma_{A_n}> \Sigma_{D_n}> \Sigma_{C_n}$, and hence
$A_n \prec_{B}D_n\prec_{B} C_n$.

Applying Lemma \ref{lem:coversecond}, we see that both $A_n$ and $D_n$ are covered by $C_n$ in the secondary Bruhat order.
So $A_n$ and $D_n$ are incomparable in the secondary Bruhat order, and the proof follows.
\end{proof}

\begin{lem}\label{lem:ngeq2k+3}
Let $k\geq 3$.
Then the Bruhat order and the secondary Bruhat order do not coincide on $\mathcal{A}(n,k)$ for all $n\geq 2k+3$.
\end{lem}

\begin{proof}
From Corollary \ref{lem:dualnkn-k} and Lemma  \ref{lem:knocoincidenk012}, we know that the Bruhat order and the secondary Bruhat order are not the same on $\mathcal{A}(k+3,k)$. Note that the Bruhat order is a refinement of the secondary Bruhat order. So there exist $(0,1)$-matrices, say $M$ and $N$, in $\mathcal{A}(k+3,k)$ such that $M\prec_{B}N$, but $M$ and $N$ are incomparable in the secondary Bruhat order.
Let
\begin{align*}
				M_{n}=\left[\begin{array}{c|c}
					M & 0 \\
					\hline 0 &  G
				\end{array}\right]\quad \text{and}\quad
N_{n}=\left[\begin{array}{c|c}
					N & 0 \\
					\hline 0 &  G
				\end{array}\right],
			\end{align*}
where $G$ is a fixed $(0,1)$-matrix in $\mathcal{A}(n-k-3,k)$.
Then $M_n\prec_B N_n$.

We claim that $M_{n}$ and $N_{n}$ are incomparable in the secondary Bruhat order on $\mathcal{A}(n,k)$.
Otherwise, we must have $M_n\prec_{\widehat{B}}N_n$ since $M_n\prec_B N_n$. Without loss of generality we may assume that $M_n$ is covered by $N_n$ in the secondary Bruhat order.
Then there exist $i,j,k,l$ such that $N_n[\{i,j\},\{k,l\}]=L_2$ and $M_n$ is the matrix obtained from $N_n$ by the $L_2\rightarrow I_2$ interchange that replaces
$N_n[\{i,j\},\{k,l\}]$ with $I_2$. It is trivial to see that $i,j,k,l\leq k+3$ so that $V_n[\{i,j\},\{k,l\}]=N_n[\{i,j\},\{k,l\}]=L_2$, and hence
$M$ can be obtained from $N$ by the $L_2\rightarrow I_2$ interchange that replaces
$N[\{i,j\},\{k,l\}]$ with $I_2$, proving $M \prec_{\widehat{B}}N $, a contradiction.
Thus, $M_{n}$ and $N_{n}$ are incomparable in the secondary Bruhat order, so the proof follows.
\end{proof}

\begin{prop}\label{lem:notequalk=4}	
The Bruhat order and the secondary Bruhat order do not coincide on $\mathcal{A}(n,4)$ for all $n$ with $n\geq 7$.
\end{prop}

\begin{proof}
It follows from Corollary \ref{lem:dualnkn-k} and Lemma \ref{lem:knocoincidenk012} that the Bruhat order and the secondary Bruhat order are not the same on $\mathcal{A}(7,4)$. Thus,
according to Lemma \ref{lem:ngeq2k+3}, we only need to consider $\mathcal{A}(n,4)$ for $n\in\{8,9,10\}$.
If $n=8$, then take
\begin{align*}
				G_1=\begin{bmatrix}
					1 & 1 & 1 & 0 \\
					0 & 0 & 0 & 1 \\
					1 & 0 & 0 & 0 \\
					0 & 1 & 1 & 1
				\end{bmatrix}, \quad G_2=\begin{bmatrix}	
					0 & 1 & 1 & 0\\
					0 & 1 & 1 & 0\\
					0 & 1 & 1 & 0\\
					1 & 0 & 0 & 1
				\end{bmatrix}\quad \text{and}\quad
			G_3=\begin{bmatrix}
				    1 & 1 & 0 & 0 \\
					0 & 0 & 1 & 1 \\
					1 & 1 & 0 & 0 \\
					0 & 0 & 1 & 1
				\end{bmatrix};
			\end{align*}
if $n=9$, then take
\begin{align*}
				G_1=\begin{bmatrix}
					1 & 1 & 1 & 0 &0\\
					0 & 0 & 0 & 0 &1 \\
					1 & 0 & 0 & 0 &0 \\
					0 & 0 & 1 & 1 &1
				\end{bmatrix}, \quad G_2=\begin{bmatrix}	
					0 & 1 & 1 & 0\\
					0 & 1 & 1 & 0\\
					0 & 1 & 1 & 0\\
					0 & 0 & 0 & 0\\
					1 & 0 & 0 & 1
				\end{bmatrix}\quad \text{and}\quad
			G_3=\begin{bmatrix}
				    0 & 0 & 0 & 1 &1\\
					0 & 1 & 0 & 1 &0 \\
					1 & 0 & 1 & 0 &0\\
					1 & 1 & 1 & 1 &0\\
					0 & 1 & 0 & 0 &1
				\end{bmatrix};
			\end{align*}
if $n=10$, then take
\begin{align*}
				G_1=\begin{bmatrix}
					1 & 1 & 1 & 0 &0 &0\\
					0 & 0 & 0 & 0 &0 &1  \\
					1 & 0 & 0 & 0 &0 &0 \\
					0 & 0 & 0 & 1 & 1 &1
				\end{bmatrix}, \quad G_2=\begin{bmatrix}	
					0 & 1 & 1 & 0\\
					0 & 1 & 1 & 0\\
					0 & 1 & 1 & 0\\
					0 & 0 & 0 & 0\\
					0 & 0 & 0 & 0\\
					1 & 0 & 0 & 1
				\end{bmatrix}\quad \text{and}\quad
			G_3=\begin{bmatrix}
				    0 & 0 & 0 & 0 & 1 &1\\
				    0 & 0 & 0 & 0 & 1 &1\\
					0 & 0 & 1 & 1 &0  &0 \\
					1 & 1 & 1 & 1 & 0 &0\\
					1 & 1 & 1 & 1 & 0 &0\\
					0 & 1 & 0 & 0 &1 &0
				\end{bmatrix}.
			\end{align*}
A simple computation shows that for any $n\in\{8,9,10\}$, the $(0,1)$-matrices
\begin{align*}
					V_n=\left[\begin{array}{c|c}
						V & G_{1}\\
						\hline G_{2}& G_{3}
					\end{array}\right]
				\end{align*}
are elements of $\mathcal{A}(n,4)$, where $V \in\{A, C, D\}$.
Therefore, by Lemma \ref{lem:xvg123}, the Bruhat order and secondary Bruhat order do not coincide on $\mathcal{A}(n,4)$.
			\end{proof}
In order to show the two Bruhat orders are not the same on $\mathcal{A}(n,k)$ for $5\leq k\leq n-3$, we need the following technical lemma.

\begin{lem}\label{lem:km-3k-22k-4k-3}
Let $k$ be a positive integer with $k\geq 5$, and let $m\in\{k,k+1,k+2\}$. If $R=(k^{m-3},(k-2)^2,(k-4)^{k-3})$,
then $R \prec R^*$.
\end{lem}
\begin{proof}
By the definition of the conjugate vector of $R$, we have
$$
R^*=((k+m-4)^{k-4},(m-1)^2,(m-3)^{2}).
$$
For notational simplicity, we write  $R =(u_{1},u_{2},\ldots,u_{k+m-4})$ and $R^{*}=(v_{1},v_{2},\ldots,v_{k})$.
It is clear that $|R|=|R^*|$ but $R\neq R^*$.
We need to show that
\begin{align}\label{eq:maineqp}
\sum_{i=1}^{s}u_i\leq \sum_{i=1}^{s}v_i
\end{align}
for all $s$ with $1\leq s\leq k+m-4$,
thereby proving  $R \prec R^*$.
Because $v_s=0$ for $k+1\leq s\leq k+m-4$, we see that $|R^*|=\sum_{i=1}^{k}v_i$.
Thus, it suffices to show Eq.~\eqref{eq:maineqp} holds for all $s$ with $1\leq s\leq k-1$.
Now the proof breaks into four steps.

{\bf Step 1.}  $1\leq s \leq k-4$. Since $k\leq m$, we have $s\leq k-4<m-3$, so that $u_s=k< k+m-4=v_s $, and hence Eq.~\eqref{eq:maineqp} follows.
	
{\bf Step 2.}  $s= k-3$. From $k\leq m$ we see that $s\leq m-3$, yielding that
\begin{align*}
\sum_{i=1}^{s}u_i=k(k-3)\quad \text{and}\quad \sum_{i=1}^{s}v_i=(k-4)(k+m-4)+(m-1)
\end{align*}
and hence
\begin{align}\label{eq:sumiq1s}
\sum_{i=1}^{s}v_i-\sum_{i=1}^{s}u_i=(k-3)m-5k+15.
\end{align}
It follows from $m\geq k\geq 5$ that
\begin{align*}
\sum_{i=1}^{s}v_i-\sum_{i=1}^{s}u_i\geq k^2-8k+15\geq 0.
\end{align*}
			
{\bf Step 3.}  $s= k-2$.
From Eq.~\eqref{eq:sumiq1s} it follows that
\begin{align*}
\sum_{i=1}^{s}v_i-\sum_{i=1}^{s}u_i=&\left[(k-3)m-5k+15\right]+(v_{k-2}-u_{k-2}).
\end{align*}
Since $v_{k-2}=m-1$ and
\begin{align*}
u_{k-2}=
\begin{cases}
k-2,&m=k,\\
k,&m=k+1\ \text{or}\ k+2,
\end{cases}
\end{align*}
we have
\begin{align}\label{eq:k-2eq}
\sum_{i=1}^{s}v_i-\sum_{i=1}^{s}u_i
=&\begin{cases}
k^2-8k+16,&m=k,\\
k^2-7k+12,&m=k+1,\\
k^2-6k+10, &m= k+2,
\end{cases}
\end{align}
which implies Eq.~\eqref{eq:maineqp} for $s=k-2$.

{\bf Step 4.}  $s= k-1$.  The proof is analogous to that of Step 3. According to Eq.~\eqref{eq:k-2eq}, together with the fact
$v_{k-1}=m-3$ and
\begin{align*}
u_{k-1}=
\begin{cases}
k-2,&m=k\ \text{or}\ k+1,\\
k,&m=k+2,
\end{cases}
\end{align*}
we obtain that
\begin{align*}
\sum_{i=1}^{s}v_i-\sum_{i=1}^{s}u_i
=&\begin{cases}
k^2-8k+15,&m=k,\\
k^2-7k+12,&m=k+1,\\
k^2-6k+9, &m= k+2,
\end{cases}
\end{align*}
from which Eq.~\eqref{eq:maineqp} follows.
Therefore, for all $s$ with $1\leq s\leq k-1$,  Eq.~\eqref{eq:maineqp} holds, so that $R\prec R^*$.
\end{proof}

\begin{proof}[Proof of Theorem \ref{thm:mainequal}]
By Corollary \ref{cor:equal} and Proposition \ref{lem:notequalk=4}, it suffices to show that the Bruhat order and the secondary Bruhat order do not coincide on $\mathcal{A}(n,k)$ for all $k$ with $3\leq k \leq n-3$.
We proceed by induction on $k$, where the base cases $k = 3$ and $k = 4$ follows from  Lemma \ref{lem:knocoincidenk012} and Proposition \ref{lem:notequalk=4}, respectively. Now assume for
positive integers $s$ less than $k$ with $3\leq s\leq n-3$, where $k\geq 5$.

If $3\leq n\leq 2k-1$, then it follows from $3\leq k \leq n-3$ that $k+3\leq n\leq 2k-1$.
Thus, by the induction hypothesis, the Bruhat order and the secondary Bruhat order do not coincide on the classes of $\mathcal{A}(n,n-k)$, and hence neither on $\mathcal{A}(n,k)$ by Corollary \ref{lem:dualnkn-k}.
If $n\geq 2k+3$, it then follows from Lemma \ref{lem:ngeq2k+3} that the two orders are also not the same on $\mathcal{A}(n,k)$.

It remains to show that the two orders are distinct on $\mathcal{A}(n,k)$ for $n\in\{2k,2k+1,2k+2\}$. Let $m=n-k$. Then $m\in\{k,k+1,k+2\}$.
By Lemmas \ref{lem:existm} and \ref{lem:km-3k-22k-4k-3}, there exists an $(n-4)$-by-$(n-4)$ $(0,1)$-matrix, say $G$, with row and column sum vectors
$$
((k-4)^{k-3},(k-2)^2,k^{m-3}).
$$
Let

\vspace{2mm}
\begin{tikzpicture}
  \draw(0,0)node{
  $V_n=\hspace{20mm}\begin{bmatrix}
					\begin{array}{llll|llllllll}
						& & & & 1 & \cdots &1 &1 &1 &0 & \cdots &0\\
						& V & & & 1&\cdots& 1 & 0&0 &0 & \cdots &0\\
						& & & & 1 &\cdots& 1 & 0&0 &0 & \cdots &0\\
						& & & & 1 & \cdots &1 &1 &1 &0 & \cdots &0\\
						\hline
						1 & 1 & 1& 1 &  & &&  & &&&\\
						\vdots & \vdots & \vdots & \vdots &  && &   & &&&\\
						1 & 1 & 1& 1 &  & &&  & &&&\\
						0 & 1 &1 &0 & &&& &&&&\\
						0 & 1 &1 &0 & &&&G &&&&\\
						0 & 0 &0 &0 & &&& &&&&\\
						\vdots & \vdots & \vdots & \vdots &  && &   & &&&\\
						0 & 0 &0 &0 & &&& &&&&\\
					\end{array}
				\end{bmatrix},$
};
\draw(0.95,3.1)node[above]{$\overbrace{\rule{16mm}{0mm}}$}(0.9,3.5)node[above]{{\scriptsize $k-3$ columns}};
\draw(4.23,3.1)node[above]{$\overbrace{\rule{16mm}{0mm}}$}(4.15,3.5)node[above]{{\scriptsize $m-3$ columns}};
\draw(-2.5,0.33)node[left]{{ {\scriptsize $k-3$ rows} $\left\{\rule{0mm}{8.1mm}\right.$ }};
\draw(-2.5,-2.33)node[left]{{ {\scriptsize $m-3$ rows} $\left\{\rule{0mm}{8.1mm}\right.$ }};
\end{tikzpicture}
\vspace{2mm}\\
where $V\in\{A,C,D\}$. Then $V_n\in \mathcal{A}(n,k)$, so the proof follows Lemma \ref{lem:xvg123}.
\end{proof}

\vspace{5mm}

\noindent
{\bf Funding}
Research partially supported by the National Natural Science Foundation of China (Grant Nos. 12071377 and 12071383) and the Educational Teaching Reform Research Project of Southwest University (Grant No. 2022JY082).
\vspace{3mm}

\noindent
{\bf Data Availability Statement}
 Data sharing not applicable to this article as no datasets were generated or
analysed during the current study.
\vspace{3mm}

\noindent
{\bf Conflict of Interests}
 The authors declare that they have no conflict of interest.
\vspace{3mm}

\noindent
{\bf  Author Contribution Statement}
Houyi Yu conceived the original idea and helped supervise this work. Tao Zhang performed the computations. Both Yu and Zhang discussed the results and contributed to the final manuscript.


\begin{thebibliography}{99}
		
\bibitem{BB05}  Bj\"orner A., Brenti F.:
        Combinatorics of Coxeter groups,
		GTM, vol. 231, Springer, New York (2005)
		
\bibitem{Bru06} Brualdi, R.A.:
                Combinatorial Matrix Classes, Encyclopedia of Mathematics and its Applications, vol. 108. Cambridge University Press, Cambridge (2006)

\bibitem{BD07} Brualdi, R.A., Deaett, L.:
               More on the Bruhat order for $(0,1)$-matrices. Linear Algebra Appl. {\bf 421}, 219--232 (2007)

\bibitem{BH04} Brualdi, R.A., Hwang S.-G.:
               A Bruhat order for the class of $(0,1)$-matrices with row sum vector $R$ and column sum vector $S$.
               Electron. J. Linear Algebra {\bf 12}, 6--16 (2004)

\bibitem{CFM13} Conflitti, A., da Fonseca, C.M., Mamede, R.:
                On the largest size of an antichain in the Bruhat order for $\mathcal{A}(2k, k)$. Order {\bf30}, 255--260 (2013)

\bibitem{CFM12} Conflitti, A., da Fonseca, C.M., Mamede, R.:
                The maximal length of a chain in the Bruhat order for a class of binary matrices. Linear Algebra Appl. {\bf 436}, 753--757 (2012)

\bibitem{FCS20} Fernandes, R., Cruz, H.F., Salom\~ao, D.:
                Classes of $(0,1)$-matrices where the Bruhat order and the secondary Bruhat order coincide. Order {\bf 37}, 207--221  (2020)

\bibitem{Gal57} Gale, D.:
                A theorem on flows in networks.
                Pacific J. Math. {\bf7}, 1073--1082 (1957)

\bibitem{Ghe14} Ghebleh M.:
                Antichains on $(0,1)$-matrices through inversions.
                Linear Algebra Appl. {\bf458}, 503--511 (2014)


\bibitem{Ry64} Ryser, H.J.:
               Combinatorial Mathematics, Carus Math. Monograph, vol. 14, Math. Assoc. of America, Providence (1964)

\bibitem{Ry57} Ryser, H.J.:
               Combinatorial properties of matrices of zeros and ones.
               Canad. J. Math. {\bf 9}, 371--377 (1957)

\bibitem{Sta12} Stanley R.P.:
                Enumerative Combinatorics, Volume I, Cambridge University Press, Cambridge (2012)
	
	\end{thebibliography}
\end{document}